\documentclass[11pt]{article}
\usepackage[a4paper,margin=3cm]{geometry}
\usepackage[utf8]{inputenc} 
\usepackage[T1]{fontenc} 
\usepackage{lmodern}
\usepackage[autolanguage]{numprint} 
\usepackage{hyperref} 
\usepackage{graphicx} 
\usepackage[all]{xy}
\usepackage{amsmath,amssymb,amsthm}

\newcommand\NN{\mathbb{N}} 
\newcommand\RR{\mathbb{R}} 
\newcommand\CC{\mathbb{C}} 
\newcommand\PP{\mathbb{P}} 

\newtheorem{theorem}{Theorem}[section]%
\newtheorem{corollary}[theorem]{Corollary}%

\newtheorem{definition}[theorem]{Definition}%
\newtheorem{remark}[theorem]{Remark}%

\author{Raphaël Butez\footnote{CEREMADE, Université Paris Dauphine, butez@ceremade.dauphine.fr}}

\title{Large deviations principle for biorthogonal ensembles and variational formulation for the Dykema-Haagerup distribution.}

\begin{document} 

\maketitle 

\begin{abstract}
This note provides a large deviations principle for a class of biorthogonal ensembles. We extend the results of Eichelsbacher, Sommerauer and Stotlz to more general type of interactions. Our result covers the case of the singular values of lower triangular random matrices with independent entries introduced by Cheliotis. In particular, we obtain as a consequence a variational formulation for the Dykema-Haagerup as it is the limit law for the singular values of lower triangular matrices with i.i.d. complex Gaussian entries.
\end{abstract}

\section{Introduction and results}
The aim of this note is to extend the work of Eichelbascher, Sommerauer and Stolz \cite{eichelsbacher2011} in order to prove a large deviations principle for a wide class of biorthgonal ensembles which include the matrix models introduced by Cheliotis in \cite{cheliotis}. The authors of \cite{eichelsbacher2011} proved a large deviations principle for wide variety of models, such as the biorthogonal Laguerre ensembles or the matrix model of Lueck, Sommers and Zirnbauer \cite{luck2006energy} for disordered bosons. Those models deal with particle systems in $\RR$ or $\CC$ with a density involving a double interaction term of type $\prod_{i<j}|x_i-x_j||x^{\theta}_i-x^{\theta}_j|$ with $\theta \in \NN^*$. Biorthogonal ensembles were introduced by Muttalib in physics in \cite{muttalib1995} and by Borodin in mathematics in \cite{borodin1998}. The recent article \cite{bloomtotik2015} develops potential theory for the model we study. Large deviations for particles systems with general repulsion have been studied in \cite{chafai} and we show that their results apply to this kind of problems.

In the article \cite{cheliotis}, Cheliotis presented a lower-triangular random matrix model for which the distribution of the singular values can be computed and form a class of biorthogonal ensembles. Later, in \cite{forrester2015}, Forrester and Wang found another matrix model for these ensembles. Large deviations principles for the empirical measures of biorthogonal ensembles enter the general framework of \cite{chafai}, but it is not clear that this model fits their technical hypotheses.

Triangular matrices are the elementary object that appear in many factorization algorithms, such at the Cholesky or the LU decomposition, so one could wonder if, starting from a random matrix, we can compute the distribution of the coefficients of it's Cholesky decomposition. Bartlett answered that question in \cite{bartlett1933theory} and proved that the entries of the Cholesky decomposition of a Wishart random matrix are independent Gaussian variables off diagonal and chi random variables on the diagonal. This result is known as the Bartlett decomposition of a Wishart matrix. Cheliotis studied the reverse problem: given a simple model of random triangular matrices $T_n$, what can we say about the eigenvalues of the eigenvalues of $T_n T_n^*$? 
	
\noindent Fix a positive integer $n \in \NN$, and two parameters $b>0$ and $\theta\geq 0$, we consider the random lower triangular matrix $$T_n = (X_{i,j})_{1 \leq i,j \leq n}$$ with independent random coefficients $X_{i,j}$ distributed according to: $$X_{i,j} \sim \begin{cases} 
\mathcal{N}_{\CC}(0,1) \text{ if } i>j ,\\ 
\frac{1}{\pi \Gamma(c_j)}e^{-|z|^2} |z|^{2(c_j-1)}d\ell_{\CC}(z) \text{ if } i=j.
\end{cases}$$
where $c_j=\theta(j-1)+b$ and  $d\ell_{\CC}$ is the Lebesgue measure on the complex plane. Note that when $\theta$ equals $0$ and $b$ equals $1$, the non-zero entries are i.i.d. complex Gaussians.

In the article \cite{cheliotis}, Cheliotis was able to compute the distribution of the ordered eigenvalues of the matrices $$S_n= T_n T_n^*.$$ He proved that, if we write $\lambda_1 \geq \dots \geq \lambda_n$ the eigenvalues of $S_n$, the distribution of the random vector $$ \Lambda_n = ( \lambda_1, \dots, \lambda_n)$$ is absolutely continuous with respect to the Lebesgue measure on $\RR^n$ with density:
\begin{equation} \label{densite theta}
\frac{1}{\prod_{j=1}^{n} j!} \frac{\theta^{-n(n-1)/2}}{\prod_{k=1}^n \Gamma(c_k)} e^{-\sum_{i=1}^{n} x_i} \prod_{j=1}^{n}x_j^{b-1}\prod_{i <j} (x_i-x_j)(x_i^\theta -x_j^\theta) 1_{x_1>\dots>x_n>0}
\end{equation}
when $\theta>0$. When $\theta=0$, the density of the distribution of $\Lambda_n$ is:
\begin{equation} \label{densite 0}
\frac{1}{\prod_{j=1}^{n} j!} \frac{1}{\prod_{k=1}^n \Gamma(c_k)} e^{-\sum_{i=1}^{n} x_i} \prod_{j=1}^{n}x_j^{b-1}\prod_{i <j} (x_i-x_j)(\log x_i-\log x_j) 1_{x_1>\dots>x_n>0}.
\end{equation}
We notice that for good choices of $\theta$ and $b$, we can recover many classical ensembles, such as the Laguerre ensembles.

In the rest of this note, we are interested in the eigenvalues of  $\frac{1}{n}S_n$. The factor $1/n$ is the proper scaling to observe a convergence of the empirical measure. We will keep the notation $\lambda_1, \dots, \lambda_n$ for the eigenvalues of $\frac{1}{n}S_n$ and we define its empirical measure:
\begin{equation*}
\mu_n = \frac{1}{n} \sum_{i=1}^{n} \delta_{\lambda_i}.
\end{equation*}

\begin{figure}[!h]
\centering
\includegraphics[scale=1]{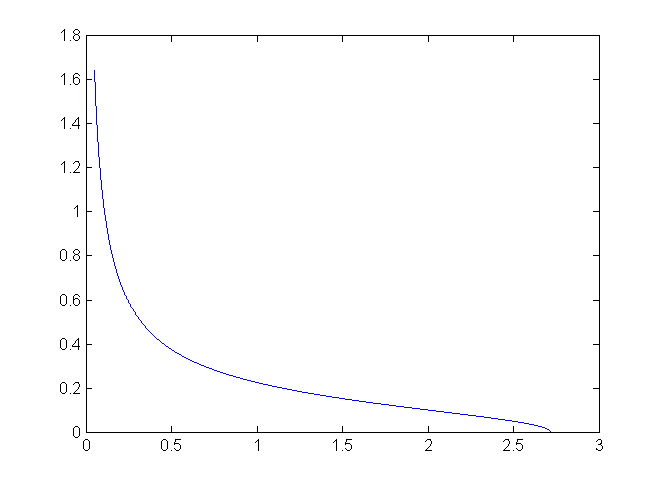}	
\caption{\label{image} Density of the Dykema-Haagerup distribution.}
\end{figure}

The special case where $\theta=0$ and $b=1$ corresponds to the case where all the coefficients are $X_{i,j}$ are independent complex random variables with variance $1$ is of particular interest.
In \cite{dykema}, using free probability theory, Dykema and Haagerup proved that $(\mu_n)_{n\in \NN^*}$ converges weakly in probability towards a deterministic measure, known as the Dykema-Haagerup distribution. In \cite{cheliotis}, the same result is proved using the moments method and path counting. The Dykema-Haagerup distribution $\mu_{DH}$ is compactly supported and absolutely continuous with respect to the Lebesgue measure on $\RR^{+*}$ with density:
\begin{equation*}
f_{DH}(x)= \frac{1}{\pi} \mathrm{Im} \left[ - \frac{1}{x W_0(x)}  \right]1_{[0,e]}
\end{equation*}
where $W_0$ is the Lambert function. $W_0$ is analytic in $\CC \setminus (-\infty, -e^{-1}]$ and can be extended to $\CC$ so that it is continuous on the upper half plane, see figure \ref{image}.

\noindent The moments of the Dykema-Haagerup distribution are given by:
\begin{equation*}
\int_{\RR} x^k d\mu_{DH}(x)= \frac{k^k}{(k+1)!}.
\end{equation*}
The Stieljes transform of $\mu_{DH}$ is defined for all $z \in \CC$ with $\mathrm{Im}(z)>$ and is given by:
\begin{equation*}
S(z)=-\frac{1}{zW_0(1/z)}=-1+e^{W_0(-1/z)}.
\end{equation*}
The R-transform of $\mu_{DH}$ is defined for all $z\in \CC$ such that $|z|<1$ and is given by:
\begin{equation*}
R(z)=-\frac{1}{(1-z)\log(1-z)}-\frac{1}{z}.
\end{equation*}

In \cite{dykema}, Dykema and Haagerup proved, using free probability tools, that the coefficients of $T_n$ are i.i.d. complex Gaussians, the empirical measures $(\mu_n)_{n\in \NN}$ converge weakly in probability towards a deterministic measure, called the Dykema-Haagerup distribution. Independently, Cheliotis proved the same result in \cite{cheliotis} using the moments method and path counting. This corresponds to the choice of $b=1$ and $\theta=0$ in our model.

In this note, the term "weak topology" corresponds to the topology associated to continuous and bounded test functions. The Bounded Lipschitz metric $d$ defined as:
\[ \forall \mu , \nu \in \mathcal{M}_1(\RR^+) \quad d(\mu,\nu)= \sup_{f} \left| \int f \mu - \int f d\nu \right| \]
where the supremum is taken over functions bounded by $1$ and $1$-Lipschitz metricizes the weak topology and makes $\mathcal{M}_1(\RR^+)$ a complete space, see \cite[Section 8.3]{bogachev2}.

\begin{definition}[Logarithmic energy]
The logarithmic energy is the functional :
\begin{equation*}
\begin{array}{cccc}
\mathcal{E} :& \mathcal{M}_1(\RR^+)& \longrightarrow &\RR \cup \{\infty\}\\ & \mu & \longmapsto & \displaystyle \iint -\log|x-y|d\mu(x)d\mu(y).
\end{array}
\end{equation*}
We also define the off-diagonal logarithmic energy
$$\begin{array}{cccc}
	\mathcal{E}_{\neq} :& \mathcal{M}_1(\RR^+)& \longrightarrow &\RR \cup \{\infty\}\\ & \mu & \longmapsto & \displaystyle \iint_{\neq} -\log|x-y|d\mu(x)d\mu(y).
\end{array}$$ where we integrate on the complement of the diagonal of $(\RR^+)^2$.
\end{definition}

\noindent We define the confining potential associated to the eigenvalue distribution \eqref{densite theta} and \eqref{densite 0}. Let $V : \RR \rightarrow \RR$ given by:
$$ V(x) = x.$$

As the eigenvalues of $\frac{1}{n}S_n$ are the eigenvalues of $S_n$ divided by $n$, we can compute the distribution of the unordered eigenvalues $(\lambda_1, \dots, \lambda_n)$. This distribution is absolutely continuous with respect to the Lebesgue measure on $\RR^n$ with density:
\begin{equation*}
\frac{1}{Z_n} e^{-n\sum_{i=1}^{n} x_i} \prod_{j=1}^{n}x_j^{b-1}\prod_{i <j} |x_i-x_j||x_i^\theta -x_j^\theta| \text{ when } \theta >0
\end{equation*}
\begin{equation*}
\frac{1}{Z_n} e^{-n\sum_{i=1}^{n} x_i} \prod_{j=1}^{n}x_j^{b-1}\prod_{i <j} |x_i-x_j||\log x_i -\log x_j|\text{ when } \theta =0
\end{equation*}
where $Z_n$ is a normalizing constant, depending on the model. Those two distributions are of the form:
\begin{equation} \label{gaz}
	\frac{1}{Z_n} e^{-n\sum_{i=1}^{n} V(x_i)} \prod_{j=1}^{n}x_j^{b-1}\prod_{i <j} |x_i-x_j||g( x_i) -g(x_j)|
\end{equation}
if we chose $g=g_{\theta}$ where:
$$g_{\theta}(x)= \begin{cases}
x^{\theta} \qquad \text{ if } \theta > 0 \\
\log x \quad \text{ else.}
\end{cases}$$
This density can be written in the form:
\begin{equation} \label{gaztheta}
\frac{1}{Z_n} \exp \left[ - n^2 \left( \frac{1}{2}\mathcal{E}_{\neq} (\mu_n) + \frac{1}{2}\mathcal{E}_{\neq} (g_*\mu_n) + \int V(x) d\mu_n(x) \right) \right]\prod_{j=1}^{n}x_j^{b-1}
\end{equation}
where $Z_n$ is a normalizing constant and where $g_{*}\mu$ is the push-forward of the measure $\mu$ by the function $g$. 

The term $\prod_{j=1}^{n}x_j^{b-1}$ will play no role in the large deviations and the same results are valid without this term. We keep this term so that the connection with the model of random matrices introduced by Cheliotis is straightforward. To recover every Laguerre ensemble, one can consider "$b=bn$", which would correspond to change the function $V$.

\begin{theorem}[Large deviations principle for $\mu_n$]\label{LDP}
Let $g$ be a $C^1$ function on $\RR^{+*}$, such that its derivative is positive. Let $V$ be a continuous function on $\RR^+$ such that there exist a constant $\beta>b$ such that we have:
\begin{equation*}
\varliminf_{x\rightarrow \infty} \frac{V(x)}{\beta \log |x|} >1 \quad \text{and} \quad \varliminf_{x\rightarrow \infty} \frac{V(x)}{\beta \log |g(x)|} >1 .
\end{equation*}
\noindent Let us define $I  : \mathcal{M}_1(\RR^+) \rightarrow \RR \cup \{\infty\}$ given by:
\[I(\mu) = 
\frac{1}{2}\mathcal{E}(\mu) + \frac{1}{2}\mathcal{E}(g_{*}\mu) + \displaystyle \int V(x)d\mu(x)\] 
The random sequence $(\mu_n)_{n\in \NN}$ satisfies a large deviations principle with speed $n^2$ in $\mathcal{M}_1(\RR)$ for the weak topology with good rate function $\tilde{I}=I_{} - \inf I_{}$. This means that for any Borel set $A \in \mathcal{M}_1(\RR^+)$ we have:
\begin{equation*}
- \inf_{\mathrm{Int}A} \tilde{I}\leq \varliminf \frac{1}{n^2}\log \PP(\mu_n \in A) \leq \limsup \frac{1}{n^2}\log \PP(\mu_n \in A) \leq -\inf_{\mathrm{Clo}A} \tilde{I}
\end{equation*}

In addition, the rate function $I_{} - \inf I_{}$ is lower semi-continuous and strictly convex on the set of the measures on which it is finite.
\end{theorem}
 \begin{remark}[Assumptions on $g$ and $V$.]
The assumptions on $g$ mean that the two interaction terms play the same role of short range repulsion, but at different scales. Our hypothesis on $g$ can be rephrased as "$g$ is locally a $C^1$-diffeomorphism of $RR^{+*}$.

The assumptions on $V$ are very standard in large deviations for Coulomb gases. They ensure that $\int e^{-V(x)}dx$ is finite and that that rate function is well defined.
 \end{remark}

In \cite{eichelsbacher2011}, Eichelsbacher, Sommerauer and Stolz proved a large deviations principle for the empirical measures $\mu_n$ when $g=g_{\theta}$ and $\theta$ is an integer and where $V$ can depend on $n$. The classical techniques to prove large deviations for the empirical measures of Coulomb gases apply here with no modification. 

The novelty of our approach is to extend the result of \cite{eichelsbacher2011} to any function $g$. Our theorem covers the original model of Muttalib from \cite{muttalib1995} with $g(x)=Argsh^2(\sqrt{x})$ which was the starting point of the study of biorthogonal ensembles. The matrix model introduced by Cheliotis corresponds to the choice of $g=g_{\theta}$ where $\theta>0$. Choosing $g(x)=\exp(x)$ gives the large deviations for the model of \cite{claeys2014}. The key of this article is the way we deal with the lower bound. Instead of inspiring from the proof of the lower bound originally given by Ben Arous and Guionnet in \cite{benarousguionnet}, we adapt the proof of Hiai and Petz from \cite{hiaipetz}. We show that the article \cite{chafai} covers a wide class of biorthogonal ensembles, which did not seem obvious. The authors of \cite{bloomtotik2015} consider a very close model as the study holomorphic functions $g$ while the density \eqref{gaz} is integrated with respect to more general measures on $\CC$ or $\RR$. Our techniques rely on the classical probabilistic approach of large deviations while they adopt a Bernstein-Markov approach.

From this result we obtain two important corollaries, which are the motivation for our study: a variational formulation for the Dykema-Haagerup distribution and the almost sure convergence of $(\mu_n)_{n \in \NN}$ towards this measure. We also state a large deviations principle for the top right particle.

\begin{corollary}[Almost sure convergence towards the minimizer]\label{convergence} 
	Let $g$ be a $C^1$ function on $\RR^{+*}$, such that its derivative is positive. 
	Let $\nu$ be the unique minimizer of the functional $I$. Then the random sequence of measures $(\mu_n)_{n \in \NN}$ converges weakly almost surely towards the deterministic measure $\nu$.
\end{corollary}

\begin{corollary}[Variational characterization of the Dykema-Haagerup Distribution] \label{Charac}
The Dykema-Haagerup distribution $\mu_{DH}$ is the unique minimizer on $\mathcal{M}_1(\RR^+)$ of the functional : $$I(\mu)= \frac{1}{2}\mathcal{E}(\mu) + \frac{1}{2} \mathcal{E}(\log_* \mu) + \int x d\mu(x)$$ which is strictly convex.
\end{corollary}

\begin{theorem}[Large deviations for the largest particle]\label{LDPtop}
Let $(x_1,\dots,x_n)$ be distributed according to \eqref{gaz} and let $x_n^*=\max_{1 \leq i \leq n} x_i$. Suppose that the hypotheses of Theorem \ref{LDP} are satisfied and assume that there exist a constant $\zeta$ such that:
\begin{equation*}
\varliminf_{n\rightarrow \infty} \frac{1}{n} \log \frac{Z^{*}_{n-1}}{Z_n} = \zeta
\end{equation*}
where $Z^{*}_{n-1}$ is the normalizing constant of the gas \eqref{gaz} with $n-1$ particles and confining potential $\frac{n}{n-1}V$.
Let $\mu_{eq}$ be the limit measure of $(\mu_n)_{n \in \NN^*}$ and let $b_eq$ be the right endpoint of its support.
The random sequence $(x_n^*)_{n \in \NN^*}$ satisfies a large deviations principle in $\RR^+$ with speed $n$ and good rate function:
\begin{equation*}
J(x) = \begin{cases}
 - \frac{1}{2} \displaystyle \int \log|x-y| + \log|g(x)-g(y)| d\mu_{eq}(y) + V(x) - \kappa&\text{ if } x \geq b_{eq}\\
 \quad \infty  &\text{ if } x < b_{eq}.\\
\end{cases}
\end{equation*}
where $\kappa$ is such that $J(b_{eq})=0$.
\end{theorem}
This theorem will not be proved in this note as the authors of \cite{credner2015} already proved this theorem for the model of \cite{eichelsbacher2011}. In the setting of \cite{credner2015}, the number of particles at step $n$ is not $n$ but $p(n)$ which makes their result more technical. One could also adapt the proof of the similar theorem from \cite{AGZ10} as the scheme of the proof is the same. First, the product structure of the density \eqref{gaz} allows us to separate the variables and integrate with respect to $x_1 < \dots < x_{n-1}$. Then, the large deviations principle for the empirical measure allows us to says that the particles $x_1 < \dots < x_{n-1}$ generate the same potential as the measure $\mu_eq$. Finally, the assumption on the normalizing constants allows us to control the error we do by changing the measure from $n$ particles to $n-1$ particles.

\begin{remark}[Large deviations for the top eigenvalue for Cheliotis' matrix model]
In the article \cite{cheliotis}, Cheliotis gives exact formulas for the normalizing constants $Z_n$ when $g=g_\theta$ and $V(x)=x$. It is straightforward to check that this model satisfies the assumptions of Theorem \ref{LDPtop}. One can also obtain another proof of the fact that for the Dykema-Haagerup model, $\lambda_{max}$ converges almost surely towards $e$.
\end{remark}

The rest of the note is devoted to the proofs of the theorems. We start by proving the large deviations principles and then we deduce the variational formula and the almost sure convergence.

\section{Proof of the large deviations principle.}
The proof of Theorem \ref{LDP} is very close to the standard proof of large deviations principle for Coulomb gases in $\RR$. Many authors proved similar results following the steps of \cite{benarousguionnet}. For general $b$ and positive integer $\theta$, theorem \ref{LDP} is a special case of the article \cite{eichelsbacher2011}. The proof is organized in several classical steps:


Not much is new in the proof that we present here, hence we will focus on what differs from the usual techniques. The parts of the proof that are omitted can be taken from \cite{AGZ10} or \cite{chafai}. The classical proof is organized as follows: \newline
Step 1: Study of the rate function;
\newline
Step 2: Exponential tightness for the non-normalized measures;
\newline
Step 3: Weak upper bound for the non-normalized measures;
\newline
Step 4: Weak lower bound for the non-normalized measures;
\newline
Step 5: Recover the full large deviations principle for the normalized measures.

We will give the fundamental inequality to prove step $1$. Then, the classical proofs of step 2 and 3 apply with no modification. We will give full details about step 4 as it is the difficult part of the proof. Once the large deviations principle is proved for the non-normalized measures, step 5 just consists in obtaining the asymptotic of the normalizing constants by applying the large deviations inequalities for the whole space of probability measure.

\subsection{Study of the rate function.}

\begin{definition} 
We set, for any non-negative $x$ and $y$:
\[ f(x,y) = -\frac{1}{2}\log|x-y| - \frac{1}{2}\log|g(x)-g(y)| + \frac{1}{2} \left[V(x) + V(y) \right]
\]
\end{definition}
\noindent Using the inequality:
\[\log|x-y| \leq \log(1+|x|)+ \log(1+|y|)\]
we obtain:
\begin{align} \label{ineq}
f(x,y) \geq  & \left( -\frac{1}{2}\log(1+|x|) - \frac{1}{2}\log(1+|g(x)|) + \frac{1}{2}V(x) \right) + \\ & \qquad  \left( -\frac{1}{2}\log(1+|y|) - \frac{1}{2}\log(1+|g(y)|) + \frac{1}{2}V(y)\right).
\end{align}

This inequality shows that the function $I$ is well defined and taxes its values in $\RR \cup \{\infty\}$. This inequality is the key to prove that $I$ is a good rate function. All the details are given in the reference book \cite[Lemma 2.6.2 p.72]{AGZ10}. 

To prove that the rate function $I$ is strictly convex where it is finite, we observe that the logarithmic energy $\mu \mapsto \mathcal{E}(\mu)$ is known to be a strictly convex function where it is finite, see \cite{AGZ10} or \cite{deift}. As the function $\mu \mapsto g_{*}\mu$ is linear, the function $\mu \mapsto \mathcal{E}(g_{*}\mu)$ is strictly convex where it is finite. The rate function $I$ is the sum of two strictly convex functions and a linear function, hence it is strictly convex on the set $\{\mu \in  \mathcal{M}_1(\RR^+) \mid I<\infty \}$.
The exponential tightness is also a consequence of the inequality \eqref{ineq}, see for instance \cite[p.77]{AGZ10}.

To prove the upper bound for non-normalized measures, the strategy of the proof is exactly the same as in \cite{chafai}. This proof applies with no modification.

\subsection{Proof of the lower bound.}
The proof of the lower bound from \cite{eichelsbacher2011} does not seem to cover the case where $g$ is not an integer power function. The classical scheme of proof from Ben Arous and Guionnet for the lower bound does not suit well for biorthogonal ensembles. We show that scheme of proof of \cite{hiaipetz} for the lower bound is more robust and allows to deal with more general types of interactions.

We want to prove that we have, for any $\sigma \in \mathcal{M}_1(\RR^+)$:
\begin{equation}
\lim_{\delta \rightarrow 0} \varliminf_{n\rightarrow \infty} \frac{1}{n^2} \log Z_n \PP(\mu_n \in B(\sigma, \delta)) \geq - I(\sigma).
\end{equation}

The classical technique consists in constructing configurations for which the density of \eqref{densite theta} is very close to $\exp(-I(\sigma))$. It is important to check that the measure of the configurations we created in $\RR^n$ does not decay too fast. Unfortunately, it is not easy to do it for general measure $\sigma$. We notice that it suffices to prove the bound for sufficiently regular measures $\sigma$. 

\subsection*{First Step: Reduction to "nice" measures.}

We will prove that for any sufficiently regular measure $\sigma$, we have:
\begin{equation}\label{eq:inf}
\inf_{G} \varliminf_{n\rightarrow \infty} \frac{1}{n^2}\log Z_n \PP(\mu_n \in G) \geq -I(\sigma)
\end{equation}
where the infimum is take over $G$ neighborhood of $\sigma$. In order to prove that this bound is sufficient to obtain the lower bound of the large deviations principle, we prove that the function $\phi : \mathcal{M}_1(\RR^+) \rightarrow \RR $ given by:
\[ \phi(\sigma) = \inf \{\varliminf_{n \rightarrow \infty} \frac{1}{n^2} \log Z_n\PP(\mu_n \in G) , G \text{ neighborhood of } \sigma\}  \] is upper semi-continuous. Let $\sigma_k \rightarrow \sigma$ in $\mathcal{M}_1(\RR^+)$. Let $G$ be a neighborhood of $\sigma$, then there exists an integer $K$ such that for all $k \geq K$, $\sigma_k \in G$. This implies that for any $k \geq K$:
\[ \inf_{G_k} \varliminf_{n \rightarrow \infty} \frac{1}{n^2} \log \PP(\mu_n \in G_k) \leq \varliminf_{n \rightarrow \infty} \frac{1}{n^2} \log \PP(\mu_n \in G)  \]
where $G_k$ are neighborhoods of $\sigma_k$. Then if we take the limit superior of this inequality and the infimum over $G$ neighborhood of $\sigma$ we obtain the upper semi-continuity of $\phi$. If we prove \eqref{eq:inf} for a dense set of measures, then for any measure $\sigma \in \mathcal{M}_1(\RR^+)$, there exist measures $\sigma_k$ such that \eqref{eq:inf} holds and $\sigma_k \rightarrow \sigma$ we get:
\begin{equation*}
\phi(\sigma) \geq \limsup_{k} \phi(\sigma_k) \geq \limsup_{k} -I(\sigma_k).
\end{equation*}
We will consider a specific sequence of measures $\sigma_k$ such that for any $k$, $\sigma_k$ is absolutely continuous with respect to the Lebesgue measure on $\RR^+$, with compact support in $\RR^{+*}$ and density bounded from above and below by positive constants and such that:
\begin{align*}
-\mathcal{E}(\sigma_k) & \geq -\mathcal{E}(\sigma) \\
\lim_{k \rightarrow \infty} \int V(x) d\sigma_k(x)  & = \int V(x) d\sigma(x).
\end{align*}
Once we obtain this sequence, we will only have to prove the lower bound for the measures satisfying the regularity conditions given above.

Let $\sigma_k= \frac{1_{1/k\leq x \leq k}}{\sigma([1/k,k])} \sigma$, then, as $f$ is bounded from below, by the monotone convergence theorem we get:
\[ \lim_{k \rightarrow \infty} \iint f(x,y) d\sigma_k(x)d\sigma_k(y) = \iint f(x,y)d\sigma(x)d\sigma(y) \]
so we can assume that $\sigma$ has compact support in $\RR^{+*}$. Now let $\phi_{\varepsilon}$ be a $\mathcal{C}^{\infty}$ probability density with support in $[0,\varepsilon]$, then we set $\sigma_{\varepsilon}= \phi_{\varepsilon} * \sigma$. The measures $\sigma_{\varepsilon}$ have compact support in $\RR^{+*}$ with continuous density and converge towards $\sigma$ as $\varepsilon$ goes to zero.

Since it is easy to check that $\int V(x) d\sigma_{\varepsilon}(x) \xrightarrow[\varepsilon \rightarrow 0]{} \int V(x) d\sigma(x)$, we only have to prove that for any $\varepsilon$
\[ - \mathcal{E}(\phi_{\varepsilon}*\sigma) \geq - \mathcal{E}(\sigma) .\]
Recall that the function $-\mathcal{E}$ is concave, so if we notice that 
\[ \phi_{\varepsilon} * \sigma = \int \phi_{\varepsilon}(y) \sigma(\cdot - y) dy \]
then, thanks to the Jensen inequality and the invariance by translation of the logarithmic energy, we obtained the desired inequality. The last thing we want for our "nice" measures is that the density is bounded from above and from below. As the density of the measures $\sigma_{\varepsilon}$ are continuous with compact support, those densities are already bounded from above. Changing $\sigma_{\varepsilon}$ to $\delta m + (1-\delta) \sigma_{\varepsilon}$ where $m$ is the uniform measure on the support of $\sigma_{\varepsilon}$ allows us to deal with measures with continuous density bounded from above and from below.

\subsection*{Second step: lower bound for "nice" measures.}
From now, $\sigma$ will be a measure with compact support $[a,b] \subset \RR^{+*}$, with density $h$ with respect to the Lebesgue measure on $\RR^+$ for which there exist a constant $C>0$ such that
\begin{equation*}
\forall x \in [a,b] \quad , \quad \frac{1}{C} \leq h(x) \leq C.
\end{equation*}
Let $a_0, \dots, a_n$ be the $\frac{1}{n}$-quantiles of $\sigma$, with $a_0=a$ and $a_n=b$. We have that for any $k$, 
\begin{equation} \label{bornes}
\frac{1}{Cn} \leq a_{k+1}-a_k \leq \frac{C}{n}.
\end{equation}
Now divide each interval $[a_{k-1},a_{k}]$ in $3$ equal parts and let $[c_k,d_k]$ be the central interval. If we set $\Delta_n = \prod_{i=1}^{n} [c_i,d_i]$, then for any $(z_1, \dots ,z_n) \in \Delta_n$, we have:
\[ d(\frac{1}{n} \sum_{i=1}^{n} \delta_{z_i}, \sigma) \leq \max_{k} |a_{k+1}-a_k| \leq \frac{C}{n} \]
where $d$ is the bounded-Lipschitz distance. We are now ready to prove the lower bound. Let $\rho_1$ be the finite measure on $\RR^+$ $x^{b-1}e^{-V(x)}dx$ and $\rho_n = \rho_1 \otimes \dots \otimes \rho_1$ the finite n-th product measure on $(\RR^+)^n$.

\begin{align*}
  & Z_n \PP(\mu_n \in B(\sigma, \delta)) \\  & \quad = \int 1_{\mu_n \in B(\sigma, \delta)} \exp \left[ - n^2 \left( \frac{1}{2}\mathcal{E}_{\neq} (\mu_n) + \frac{1}{2}\mathcal{E}_{\neq} (g_{*}\mu_n) + 																										\frac{n-1}{n}\int V(x) d\mu_n(x) \right) \right] d\rho_n(x)\\
  & \quad \geq \int 1_{\Delta_n} \exp \left[ - n^2 \left( \frac{1}{2}\mathcal{E}_{\neq} (\mu_n) + \frac{1}{2}\mathcal{E}_{\neq} (g_{*}\mu_n) + \frac{n-1}{n}\int V(x) d\mu_n(x) \right) \right] d\rho_n(x) \\ 
  & \quad \geq \exp \left( -n^2\left[\frac{n-1}{n^2}\sum_{k=1}^{n} \max_{ [c_{i},d_{i}]}V(x) \right] \right) \times  \exp \left( -n^2\left[ - \frac{1}{n^2} \sum_{i < j} \min_{[c_{i},d_{i}] \times [c_{j},d_{j}]} \log|x-y| \right] \right) \times \\ & \qquad  \exp \left( -n^2\left[- \frac{1}{n^2} \sum_{i < j} \min_{[c_{i},d_{i}] \times [c_{j},d_{j}]} \log|g(x)-g(y)|   \right] \right)\int 1_{\Delta_n} d\rho_n(x).
\end{align*}
We notice that:
\[\frac{1}{n^2} \log  \int 1_{\Delta_n} d\rho_n(x) \xrightarrow[n \rightarrow \infty]{} 0. \]
Hence, to obtain the lower bound, it is sufficient to prove that we have:
\begin{equation}\label{V}
\lim\limits_{n\rightarrow \infty} \frac{1}{n} \sum_{k=1}^{n} \max_{[c_{i},d_{i}]}V = \int V(x) d\sigma(x),
\end{equation} 
and, using the fact that the functions logarithm and $g$ are increasing:
\begin{equation} \label{minor}
\varliminf_{n \rightarrow \infty} \frac{1}{n^2} \sum_{i < j} -\log(d_j-c_i) \geq \frac{1}{2}\iint -\log|x-y|d\sigma(x)d\sigma(y)=\frac{1}{2} \mathcal{E}(\sigma)
\end{equation}
and also:
\begin{equation} \label{minorg}
\varliminf_{n \rightarrow \infty} \frac{1}{n^2} \sum_{i < j} \log(g( d_j) - g( c_i)) \geq  \frac{1}{2}\iint \log|g( x)-g( y)|d\sigma(x)d\sigma(y)= \frac{1}{2}\mathcal{E}(g_{*} \sigma).
\end{equation}

If we admit temporarily the inequalities \eqref{V}, \eqref{minor} and \eqref{minorg}, the proof of the lower bound for regular measures is completed. The last step will consist in proving those three inequalities.
\subsection*{Last step: Proof of the inequalities.}

First, \eqref{V} is easy to check as we approximate a continuous integrable function on $[a,b]$ by simple functions.

We now prove \eqref{minor} following the proof of \cite{hiaipetz}. We admit temporarily that there exist a constant $A>0$ such that for $i<j$:
\begin{equation}\label{eq:borne}
A (d_j-c_i) \geq  (a_j-a_{i-1})
\end{equation}
and also that:
\begin{equation}\label{eq:proportion}
\lim\limits_{n \rightarrow \infty} \frac{2}{n^2}\# \{ i < j \mid \frac{(a_j-a_{i-1})}{(d_j-c_i)}\leq 1+\varepsilon \}=1. 
\end{equation} 
We postpone the proof of the inequalities \eqref{eq:borne} and \eqref{eq:proportion} to prove \eqref{minor}. We call:
$$B_n= \mathcal{E}(\sigma) - \frac{2}{n^2}\sum_{i\neq j}  \left(\min_{[c_i,d_i]\times[c_j,d_j]}\log|x-y|\right) $$
and we want to prove that:
\begin{equation*}
\varliminf_{n\rightarrow \infty} B_n \leq 0
\end{equation*}
Since
\begin{equation*}
\iint \log |z-w|d\sigma(z) d\sigma(w) \leq  \frac{2}{n^2}\sum_{i < j} \log|a_j-a_{i-1}| + \frac{1}{n^2} \sum_{i=1}^{n} \log |a_i - a_{i-1}|
\end{equation*} 
then for every $\varepsilon>0$ we have :
\begin{align*}
B_n   \leq & \frac{2}{n^2}\sum_{i< j} \log |a_j - a_{i-1}|  - \frac{2}{n^2}\sum_{i< j} \log |d_j - c_{i}|+ \frac{1}{n^2} \sum_{i=1}^{n} \log |a_i - a_{i-1}|\\  
\leq & \frac{2}{n^2}\# \{ i < j \mid \frac{(a_j-a_{i-1})}{(d_j-c_i)}\leq 1+\varepsilon \} \log(1+\varepsilon)  \\  & \quad  +  \frac{1}{n^2}\left[1-\frac{2}{n^2}\# \{ i < j \mid \frac{(a_j-a_{i-1})}{(d_j-c_i)}\leq 1+\varepsilon \}\right] \log A +\frac{1}{n^2} \sum_{i=1}^{n} \log |a_i - a_{i-1}|.
\end{align*}
Then we take the limit superior in both sides, and the limit when $\varepsilon \rightarrow 0$
\begin{align*}
\mathcal{E}(\sigma) - \varliminf_{n\to\infty} \frac{2}{n^2}\sum_{i< j} \log |d_j-c_{i}| \leq 0
\end{align*}
which proves \eqref{minor}. 

We prove now inequality \eqref{eq:borne}. From inequality \eqref{bornes}, we get for any $k>0$:
\begin{equation*}
\frac{a_{i+k}-a_{i-1}}{d_{i+k}-c_{i}} \leq \frac{(k+1)C/n}{(k+2/3)/Cn}.
\end{equation*}
We deduce from this inequality that the left part of the inequality is bounded by a constant independent of $k$ and $n$, which proves \eqref{eq:borne}. 
In order to prove \eqref{eq:proportion}, we start from:
\begin{equation*}
\frac{a_{i+k}-a_{i-1}}{d_{i+k}-c_{i}}= 1+ \frac{a_{i+k}-d_{i+k}}{d_{i+k}-c_{i}} +\frac{c_{i}-a_{i-1}}{d_{i+k}-c_{i}}.
\end{equation*}
If we show that: 
\[ \frac{a_{i+k}-d_{i+k}}{d_{i+k}-c_{i}} \quad \text{and} \quad \frac{c_{i}-a_{i-1}}{d_{i+k}-c_{i}} \]
can be made as small as desired when $k$ is bigger than a certain constant independent of $n$, then \eqref{eq:proportion} will be proved. 
Using \eqref{bornes} we get:
\begin{equation*}
\frac{a_{i+k}-d_{i+k}}{d_{i+k}-c_{i-1}} \leq \frac{C/3n}{k/Cn}
\end{equation*}
and
\begin{equation*}
\frac{c_{i}-a_{i-1}}{d_{i+k}-c_{i}} \leq \frac{C/3n}{k/Cn}.
\end{equation*}
Those two terms can be made as small as desired is $k$ is sufficiently large, independently of $n$, which proves \eqref{eq:proportion}.

The proof of the inequality \eqref{minorg} mimics the proof of inequality \eqref{minor}. Like in the previous case, it is sufficient to find a constant $A'$ such that for any $i<j$:
\begin{equation} \label{eq:borneg}
A' (g(d_j)-g(c_{i-1})) \geq g(a_j)-g(a_{i-1})
\end{equation}
and to prove that: 
\begin{equation} \label{eq:proportiong}
\lim\limits_{n \rightarrow \infty} \frac{2}{n^2}\# \{ i < j \mid \frac{g(a_j)-g(a_{i-1})}{g(d_j)-g(c_{i}))}\leq 1+\varepsilon \}=1. 
\end{equation} 
As the support of $\sigma$ is a compact included in $\RR^{+*}$, there exist two constants $m$ and $M$ such that for all $x  \in [a,b]$:
\[m \leq g'(x)\leq M. \] 
The inequality \eqref{eq:borneg} is a consequence of \eqref{eq:borne}, using the mean value theorem for $g$ and the fact that its derivative is bounded from above and from below.
\noindent The inequality \eqref{eq:proportiong}  is equivalent to prove that the quantities
\[ \frac{g(a_{i+k})-g(d_{i+k})}{g(d_{i+k})-g(c_{i})} \quad \text{and} \quad \frac{g(c_i)-g(a_i)}{g(d_{i+k})-g(c_{i})}\]
are as small when $k$ is large enough. Using the mean value theorem we get:
\begin{equation*}
\frac{g(a_{i+k})-g(d_{i+k})}{g(d_{i+k})-g(c_{i})} \leq \frac{M}{m} \frac{a_{i+k}-a_{i-1}}{d_{i+k}-c_{i}} \leq \frac{M}{m}\frac{(k+1)C/n}{(k+2/3)/Cn}.
\end{equation*}
The other term is treated in the same way.
\noindent Now that we have proved \eqref{eq:borneg} and \eqref{eq:proportiong}, the proof of \eqref{minorg} is the exactly the same as the proof of \eqref{minor}.

\section{Proof of Corollary \ref{convergence} and Corollary \ref{Charac}.}
As the function $I$ is lower semi-continuous and strictly convex, it has a unique minimizer, called $\mu_{b,/theta}$. Consider the sets:
\[A_{\varepsilon} = \mathcal{M}_1(\RR^+) \setminus B(\nu,\varepsilon) .\]
As $I$ is lower semi-continuous, $\inf\{ I_(\mu) , \mu \in A_{\varepsilon} \} >0$, then, thanks to the Borel-Cantelli lemma, we get:
\[ d(\mu_n,\nu) \xrightarrow[n \rightarrow \infty]{} 0.\]

As we already know that when $b=1$ and $\theta=0$ the random sequence $(\mu_n)_{n\in \NN^*}$ converges weakly in probability towards the Dykema-Haagerup distribution $\mu_{DH}$. We also know from Corollary \ref{convergence} that $(\mu_n)_{n\in \NN^*}$ converges almost surely weakly towards the minimizer of $I$. Hence we obtain the following characterization of $\mu_{DH}$:

\begin{equation*}
\mu_{DH}= \inf \{ \frac{1}{2} \mathcal{E}(\mu) + \frac{1}{2}\mathcal{E}(\log_* \mu)  +\int x\mu(x) , \mu \in \mathcal{M}_1(\RR^+)  \}.
\end{equation*}

\section{Perspectives.}
We can extend to any finite number of interactions of type: $$\prod_{i<j}|f_1(x_i)-f_1(x_j)|^{\beta_1}\prod_{i<j}|f_2(x_i)-f_2(x_j)|^{\beta_2} \dots \prod_{i<j}|f_p(x_i)-f_p(x_j)|^{\beta_p}$$
where each of the $f_k$ is locally a $C^1$ diffeomorphism and the $\beta_k$ are positive numbers. Large deviations will be valid if the confining potential $V$ dominates all the functions $f_k$ at the same time at infinity. The proof of this result would be similar to the proof of Theorem \ref{LDP}.

The result of this note can be extended in any dimension if we make additional assumptions on the function $g$. One could assume that $g$ is continuously differentiable and that on any compact $K$, there exist a constant $m_K$ such that for any $x,y \in K$:
\begin{equation*}
\| g(x)-g(y) \| \geq m_K \| x-y \|.
\end{equation*}
This condition is equivalent to $g$ being  locally a $C^1$ diffeomorphism. The article \cite{bloomtotik2015} covers the complex case.

The model studied by Götze and Vencker in \cite{gotzevenker} is not covered by this note, as they deal with a double interaction term of the type $\prod_{i<j}|x_i-x_j|^2\phi(x_i-x_j)$. This is really the combination of two different interactions whereas our model deals with the usual logarithmic interaction at two different scales. As this model is covered by the study \cite{chafai}, one could try to find the optimal conditions of $\phi$ so that a large deviations principle is valid.

We would like to thank Dimitris Cheliotis whose work \cite{cheliotis} is the starting point of this study. 
\bibliographystyle{alpha} 
\bibliography{biblio} 
\end{document}